\documentclass[12pt]{article}

\usepackage{url}
\usepackage{fancyhdr}
\usepackage{setspace}  
\usepackage{mathptmx}
\usepackage{latexsym}
\usepackage{amssymb}
\usepackage{amsthm}
\usepackage{amsmath}
\usepackage{authblk}
\newtheorem{theorem}{Theorem}
\newtheorem{lemma}{Lemma}
\newtheorem{corollary}{Corollary}

\begin{document}

\title{All face 2-colorable d-angulations \\are Gr\"unbaum colorable}

\author{\Large{Serge Lawrencenko}}

\author{\Large{Abdulkarim M. Magomedov}}
\affil[]{}
\renewcommand\Authands{}
\maketitle

\begin{abstract}
A $d$-angulation of a surface is an embedding of a 3-connected graph on that surface that divides it into $d$-gonal faces. A $d$-angulation is said to be Gr\"unbaum colorable if its edges can be $d$-colored so that every face uses all $d$  colors. Up to now, the concept of Gr\"unbaum coloring has been related only to triangulations ($d = 3$),  but in this note, this concept is generalized for an arbitrary face size $d \geqslant 3$. It is shown that the face 2-colorability of a {$d\mbox{-}$angulation} $P$  implies the Gr\"unbaum colorability of $P$.  Some wide classes of triangulations have turned out to be face 2-colorable.  
\end{abstract}


\hphantom{10 pt} \hphantom{10 pt} \hphantom{1 pt} \hphantom{1 pt} \hphantom{10 pt} 05B45, 52C20, 57M20, 57M15 (Secondary).

\section{Terminology and notation}
In this note we consider only {\it simple} graphs---that is, graphs without loops or parallel edges, and we embed the graphs only on {\it closed} surfaces---that is, surfaces without boundaries, such as a sphere. We mainly follow the standard terminology and notation of graph theory (\cite{H}). 

A {\it $d$-angulation} $P$  of a surface means a {\it $d$-gonal embedding} of a 3-connected graph $G=G(P)$  on that surface---that is, an embedding each face of which is bounded by a simple cycle of $G$ with fixed length 
$d \geqslant 3$.  Combinatorially, $P$  is defined by the triple of sets $V(P), E(P)$, and $F(P)$  of vertices, edges, and faces, respectively. The {\it dual graph} $G^*(P)$  is defined to be the graph the vertex sets of which corresponds to $F(P)$ and in which two vertices are adjacent if and only if the corresponding faces of $P$  are adjacent. Notice that $G^*(P)$ is {\it $d$-regular}---that is, the degree of each vertex is equal to $d$. 

In this note we only consider $d$-angulations $P$  whose dual graphs are simple graphs, and therefore we will suppose that $G(P)$  is 3-connected, which ensures the simplicity of $G^*(P)$  in the following two important cases: 
\par\medskip
\noindent (i)~~$d = 3,$ 

\noindent (ii)~the carrier surface is a sphere. 
\par\medskip
In case ii, by Steinitz's Theorem, every 3-connected planar graph $G$ is the $1\mbox{-}$skeleton of a convex polytope (in 3-space) with boundary complex $P$,  and the dual graph $G^*(P)$ appears to be the 1-skeleton of the dual polytope with boundary complex $P^*$. 

A {\it vertex}, {\it edge}, or {\it face $k$-coloring} of a $d$-angulation $P$ is a surjection of the set $V(P)$,  $E(P)$, or $F(P)$  onto a set of $k$ distinct colors such that the images of adjacent vertices, edges, or faces are different, respectively. Especially, any edge 3-coloring of a 3-regular graph is called a {\it Tait coloring}. 
The {\it vertex}, {\it edge}, and {\it face chromatic numbers} of $P$ are defined to be the smallest values of $k$ possible to obtain corresponding $k$-colorings, and are denoted by $\chi(P)$, $\chi^\prime(P)$,  and $\chi^{\prime\prime}(P)$, respectively. The numbers $\chi(P)$ and $\chi^\prime(P)$  are also called the vertex and edge chromatic numbers of the graph $G(P)$  itself and denoted as $\chi(G(P))$  and $\chi^\prime(G(P))$,  respectively.

Clearly, any face $k$-coloring of an arbitrary $d$-angulation $P$  corresponds to some vertex $k$-coloring of the dual graph $G^*(P)$, and conversely, whence 
$$\chi^{\prime\prime}(P)=\chi(G^*(P)).$$  
Interestingly, since $G^*(P)$  is $d$-regular, it follows that 
$${{\chi(G^*(P)) \in \{2, 3, \ldots, d, d+1 \}}},$$ 
and that 
$$\chi ^{\prime} (G^* (P)) \in \{3, \ldots, d, d+1 \}$$ 
by Vizing's Theorem \cite{V}.

\section{Gr\"unbaum colorings}
 A {\it Gr\"unbaum coloring} is a coloring of the edges of a $d$-angulation $P$  with $d$  colors such that  for each face $f$  all $d$  colors occur at the edges incident to $f$.  Up to now, the concept of Gr\"unbaum coloring has been related only to triangulations---that is, the case $d = 3$, but in this note we will generalize this concept for an arbitrary face size $d \geqslant 3$.  If $T$  is a triangulation, then $\chi^{\prime}(G^*(T)) \in \{ 3,4\}$  by Vizing's Theorem. The equality $\chi^{\prime}(G^*(T))=3$  means that $G^*(T)$  is Tait colorable, 
which, in the dual form, means that $T$  is Gr\"unbaum colorable. 
\par\medskip    
\par\medskip 
\noindent {\bf Conjecture 1~(Gr\"unbaum \cite{G}, 1969).}~Every triangulation $T$  of an orientable surface is Gr\"unbaum colorable---that is, $\chi^{\prime}(G^*(T))=3$.
\par\medskip    
\par\medskip 

Conjecture 1 stood for 40 years, until Kochol \cite{Koch} constructed infinite families of counterexamples on orientable surfaces with genus $g$ for all $g \geqslant 5.$  Here we put forward another conjecture about triangulations by strengthening the vacuous restriction $\chi^{\prime\prime}(T)\leqslant 4$  (which obviously holds for any $T$) to the restriction $\chi^{\prime\prime}(T)\leqslant~3$, but without restricting the surface's orientability class:
\par\medskip    
\par\medskip   

\noindent {\bf Conjecture 2.}~Every triangulation $T$  of (a) an orientable or (b) nonorientable surface with $\chi^{\prime\prime}(T)\leqslant 3$  is Gr\"unbaum colorable.
\par\medskip    
\par\medskip  

In Section 3 we establish (Theorem 2) that for Gr\"unbaum colorability of a $d$-angulation $P$  (that is, for the equality $\chi^{\prime}(G^*(P))=d$  to hold), it suffices that $\chi^{\prime\prime}(P) = 2$,  without the orientability restriction. In Sections 4 and 5 we establish the face 2-colorability in some known, and quite wide, classes of triangulations.

\section{Key Theorem}

Let $P$  be a $d$-angulation of an orientable or nonorientable surface (whose dual graph is a simple graph). Since the dual graph $G^*(P)$  is $d$-regular, the following lemma is obvious.

\begin{lemma}
In order for the equality $\chi^{\prime}(G^*(P))=d$  to hold, it is necessary and sufficient that the graph $G^*(P)$  be 1-factorable---that is, be the sum of $d$  one-factors.
\end{lemma}

In a classical article, K\"onig \cite{Ko} (also see \cite{LP}) proved that each bipartite $d\mbox{-}$regu-lar graph expands to the sum of $d$  one-factors. Since a graph is bipartite if and only if it is vertex 2-colorable, we get the following reformulation of K\"onig's Theorem:

\begin{theorem} [\bf K\"onig]
If $\chi(G^*(P))=2$, then $G^*(P)$  is 1-factorable.
\end{theorem}
   
By a combination of Lemma 1 and Theorem 1, we obtain our key theorem which states that each face 2-colorable $d$-angulation of an orientable or nonorientable surface is Gr\"unbaum colorable: 

\begin{theorem} [\bf Key Theorem]
If $\chi(G^*(P))=2$,  then $\chi^{\prime}(G^*(P))=d$. \\
\indent\indent\indent{\rm {Dual formulation:}} If $\chi^{\prime\prime}(P)=2$,  then $P$  is Gr\"unbaum colorable.
\end{theorem}

As a particular case of Theorem~2, when $d=3$,  we can state that Conjecture~1 certainly holds for all face 2-colorable triangulations of orientable and nonorientable surfaces. 
Notice that in Theorem 2 the face chromaticity condition is only minimally strengthened in comparison to that in Conjecture 2.  

Conjecture 1 in full generality is obviously false if extended to the nonorientable case. The best known counterexample is provided by the minimal triangulation $T_{\min}$  of the projective plane by the complete 6-graph $G=K_6$. In this case, $G^*(T_{\min})$  turns out to be the Petersen Graph \cite{P} which cannot be decomposed into the sum of three 1-factors (see \cite{H}, \cite{LP}, \cite{P}) and by Lemma 1 has edge chromatic number at least 4. (This number is in fact equal to 4, which, by the way, easily implies the non-Hamiltonicity of the Petersen Graph---these are excellent creative exercises for a college course on Discrete Mathematics!)

\section{Triangulations by complete graphs}

In this section, we establish the existence of Gr\"unbaum colorable triangulations on orientable and nonorientable surfaces by complete graphs $K_n$  for at least half of the residue classes in the spectrum of possible values of $n$.
     
We begin with the orientable case, in which there exists a triangulation by $K_n$  if and only if $n\equiv 0, 3, 4$ or 7 (mod 12); see \cite{R}. Grannell, Griggs, and \v{S}ir\'a\v{n} \cite{GGS}  noticed that, when $n\equiv 0$  or 4 (mod 12), such triangulations are not face 2-colorable because for face 2-colorability it is necessary that all vertex degrees should be even, that is, $n$  should be odd. Furthermore, they established that the orientable triangulations constructed by Ringel \cite{R}  for all $n\equiv 3$ (mod 12)  are face 2-colorable; however, since $K_3$ is not 3-connected, we have to enforce $n \ne 3$ (see Section 1).  Finally, they established that there exists a face 2-colorable triangulation for each $n\equiv 7$ (mod 12) among the orientable triangulations constructed by Youngs \cite{Y}. We summarize these results in the following theorem. 

\begin{theorem}  [\bf \bf Ringel \cite{R}; Youngs \cite{Y}; Grannell, Griggs, \v{S}ir\'a\v{n} \cite{GGS}]
There exists a face $2$-colorable triangulation of an orientable surface by the complete graph $K_n$  if and only if $n\equiv 3$  or $7$ {\rm (mod 12)},  $n \ne 3$.
\end{theorem}

If one triangulation is face 2-colorable and the other is not, the two triangulations are certainly  \textit{nonisomorphic}---that is, there is no bijection between their vertex sets that extends to a homeomorphism between the surfaces carrying the triangulations. 

Historically, the first examples of pairs of nonisomorphic orientable triangulations with the same complete graph were constructed \cite{Y} in 1970. In those examples, the non-isomorphism follows immediately from the fact that one of the triangulations is face 2-colorable while the other is not; see review \cite{GG}. After a quarter of a century, in \cite{LNW}, there was constructed an example of {\it more than two} nonisomorphic orientable triangulations with the same complete graph, namely: there were constructed {\it three} such triangulations, only one of which is face 2-colorable. In 2000, it was shown \cite{BGGS} (also see \cite{GG}) that the number of nonisomorphic orientable triangulations with graph $K_n$  actually grows very rapidly as $n \to \infty$  even within the class of face 2-colorable triangulations; for instance, when $n\equiv 7$ or $19$ (mod $36$), that number is at least $2^{n^2/54-o(n^2)}$. 

The following corollary can be proved by a combination of Theorems 3 and 2.

\begin{corollary}
For each $n\equiv 3$ or $7$ {\rm (mod 12)},  $n \ne 3$,  there exists a Gr\"unbaum colorable orientable triangulation by the complete graph $K_n$.
\end{corollary}
   
To turn to the nonorientable case, recall \cite{R} that $K_n$  triangulates a nonorientable surface if and only if $n\equiv 0$ or $1$ (mod 3), $n \geqslant 6$  and $n \ne 7$.  Thus, under these conditions, $n$  is odd if and only if $n\equiv 1$  or $3$ (mod $6$), $n\geqslant 9$.  For all these values of $n$,  face $2$-colorable triangulations of the corresponding nonorientable surface by the graph $K_n$ are constructed in \cite{R} and \cite{GKo}.

\begin{theorem} [\bf Ringel \cite{R}; Grannell, Korzhik \cite{GKo}]
There exists a face $2$-colorable triangulation of a nonorientable surface by the complete graph $K_n$  if and only if $n\equiv 1$  or $3$ {\rm (mod 6)}, $n\geqslant 9$.
\end{theorem}

The following corollary can be proved by a combination of Theorems 4 and 2.
  
\begin{corollary}
For each  $n\equiv 1$ or $3$ {\rm (mod 6)}, $n\geqslant 9$,  there exists a Gr\"unbaum colorable nonorientable triangulation by the complete graph $K_n$.
\end{corollary}

Theorems 3 and 4 guarantee that we have not missed any face 2-colorable triangulations when applying Theorem 2 for obtaining Corollaries 1 and 2 (respectively).

\section{Triangulations by tripartite graphs}
Firstly, we notice that the existence of an orientable triangulation by the complete tripartite graph $K_{n,n,n}$   was established by Ringel and Youngs \cite{RY} for each $n$.  Secondly, the face 2-colorability of each such triangulation was established by Grannell, Griggs, and Knor \cite{GGK} (also see \cite{GG}). A combination of these two results with Theorem 2 leads to the following statement: {\it  For each $n \geqslant 2$,  all triangulations of the corresponding orientable surface by the complete tripartite graph $K_{n,n,n}$  are Gr\"unbaum colorable.} However, as observed by Archdeacon \cite{A}, it is very easy to prove this fact directly, even without using the completeness or orientability conditions: if the vertex parts are $A$, $B$, $C$, then color the edges between $A$  and $B$  red, those between $B$  and $C$  blue, and those between $A$  and $C$  green, and we are done!

At first sight, the statement of the preceding paragraph may seem to be subjectless; however, as shown in \cite{GGKS} (see also \cite{GG}) in the case $n$  is prime, there exist at least $(n-2)!/(6n)$  nonisomorphic orientable triangulations by $K_{n,n,n}$.  Furthermore, \cite{GKn}  provides improved bounds on the number of such triangulations; for instance, when $n\equiv 6$ or $30$ (mod $36)$,  there exist at least $n^{n^2/144-o(n^2)}$  nonisomorphic orientable triangulations by $K_{n,n,n}$.

\section*{Acknowledgment}
The first author wishes to thank Dan Archdeacon for helpful discussions.

\bibliographystyle{model1-num-names}
\bibliography{<your-bib-database>}

\par\medskip 
\par\medskip 
\par\medskip 

\noindent \textsc{S. Lawrencenko} \\
\noindent Russian State University of Tourism and Service,\\
\noindent Lyubertsy, Moscow Region, Russia,\\
\noindent e-mail: \url{lawrencenko@hotmail.com}\\

\noindent \textsc{A.\,M. Magomedov} \\
\noindent Dagestan State University, Makhachkala, Dagestan, Russia,\\
\noindent e-mail: \url{magomedtagir1@yandex.ru}

\end{document}